\documentclass[12pt]{elsarticle}

\renewcommand{\wr}{\mathbin{\rm wr}}
\newcommand\Aut{\mathop{\mathrm{Aut}}\nolimits}
\newcommand\Prob{\mathop{\mathrm{Prob}}\nolimits}
\def\Out{\mathop{\mathrm{Out}}\nolimits}
\def\PSL{\mathop{\mathrm{PSL}}\nolimits}
\def\PGL{\mathop{\mathrm{PGL}}\nolimits}
\def\PGaL{\mathop{\mathrm{P}\Gamma\mathrm{L}}\nolimits}
\newcommand\eee{\mathrm e}

\def\half{{1/2}}
\newtheorem{thm}{Theorem}[section]
\newtheorem{cor}[thm]{Corollary}

\newtheorem{lem}[thm]{Lemma}
\newproof{pf}{Proof}
\newdefinition{rk}{Remark}

\begin{document}

\begin{frontmatter}

\title{Most primitive groups are full
automorphism groups of edge-transitive hypergraphs} 
\author{L\'aszl\'o Babai}
\address{Department of Computer Science, University of Chicago, 1100 E 58th St.,
Chicago~IL~60637, USA}
\author{Peter J. Cameron\corref{cor}}
\ead{pjc20@st-andrews.ac.uk}
\address{School of Mathematics and Statistics, University of St Andrews, North
Haugh, St~Andrews, Fife KY16 9SS, UK}
\cortext[cor]{Corresponding author}

\begin{abstract}
We prove that, for a primitive permutation group $G$ acting on a set $X$ of
size $n$, other than the alternating group, the probability that
$\Aut(X,Y^G) = G$ for a random subset $Y$ of $X$, tends to $1$ as $n
\rightarrow \infty$. So the property of the title holds for all
primitive groups except the alternating groups and finitely many others.
This answers a question of M.~H.~Klin. Moreover, we give an upper bound
$n^{1/2+\epsilon}$ for the minimum size of the edges in such a hypergraph.
This is essentially best possible.
\end{abstract}


\end{frontmatter}

\centerline{\em Dedicated to the memory of \'Akos Seress}

\section{Introduction}

	It is well known that, although every abstract group is the
full automorphism group of a graph, not every permutation group is.
Moreover, the alternating group is not the automorphism group of any
family of sets, or even the intersection of automorphism groups of
families of sets; for any set system admitting the alternating group
admits the symmetric group. It is our purpose here to show that, at
least for primitive groups, there are only finitely many other
exceptions. Moreover, we can assume that $G$ acts transitively on the
sets of the system, and that the size of the sets is not too large.

\begin{thm}  
For any $\epsilon>0$, there is a finite list $\mathcal{L}$ of primitive
permutation groups 
such that the following holds for every $k$ satisfying
$n^{1/2+\epsilon} \le k \le n/2$.
If $G$ is a primitive
permutation group of degree $n$, which is not the alternating group
and not in the list $\mathcal{L}$, then there is a
$k$-uniform 
hypergraph
$(X,\mathcal{B})$
such that
$\Aut(X,\mathcal{B})=G$ and $G$ acts transitively on $\mathcal{B}$.
\label{best}\end{thm}

\begin{cor}
Let $G$ be a primitive group on $X$, not the alternating group and not
one of a finite list of exceptions. Then there is a uniform hypergraph
$(X,\mathcal{B})$ such that $\Aut(X,\mathcal{B}) = G$ and $G$ acts
transitively on $\mathcal{B}$.
\label{cor1}\end{cor}

If we do not restrict the edge size, then 
the $G$-orbits of almost all subsets of the vertex set define 
edge-transitive hypergraphs that  
have automorphism group no larger than $G$:

\begin{thm}
Let $|X| = n$, and let $G$ be a primitive permutation group on $X$ but
not the alternating group. If \ $Y$ is a random subset of $X$ and $Y^G$
the set of $G$-translates of $Y$ then
\[
\Prob(\Aut(X,Y^G)>G) < \exp(-n^{\half + o(1)}).
\]
\label{all}\end{thm}

\begin{rk}
We give first the proof of Theorem~\ref{all}, since
Theorem~\ref{best} uses similar arguments but needs more refined
estimates.
\end{rk}

\begin{rk}
We have not attempted to determine the ``finite list of
exceptions'' in Corollary~\ref{cor1}. Note that any set-transitive group is
an exception. There are just four of these apart from $S_n$ and $A_n$,
\textit{viz.}\ the Frobenius group of order $20$ ($n=5$), $\PGL(2,5)$ ($n =6$),
$\PGL(2,8)$ and $\PGaL(2,8)$ ($n=9$). Another exception is the Frobenius
group of order $21$ ($n=7$); any orbit (or union of orbits) of $G$ on
3-sets (or on 4-sets) admits one the three minimal overgroups of
$G$ (the Frobenius group of order 42 and one of two copies of $\PGL(3,2)$).

A very similar situation arose in connection with the main theorem of
\cite{cns}, where it was shown that every primitive group apart from symmetric
and alternating groups and a finite list has a regular orbit on the power
set of its domain. The finite list was computed by Seress~\cite{seress}.
His methods were our inspiration to complete the work reported here.
\end{rk}

\begin{rk}
What is the least size of edges in a hypergraph $(X,\mathcal{B})$
with $\Aut(X,\mathcal{B}) = G$? We cannot get by with edges of fixed size.
Consider, for example, the alternating group $G$ of degree $m$ in its
induced action on 2-sets. A $k$-subset $Y$ of $X$ is the edge set of a
graph with $m$ vertices and $k$ edges. If $\mathcal{B}$ is a collection of
$k$-sets with $\Aut(X,\mathcal{B})=G$, then some set $Y\in\mathcal{B}$ does not
admit any odd permutation in $S_m$; so it has at most one fixed point, and at
least $(m-1)/2=\Omega(\sqrt{n})$ edges.
This shows that Theorem~\ref{best} is best possible, apart from the value
of $\epsilon$ in the exponent.
\end{rk}

Our proof uses
the following result \cite{fpgfsg}:
\begin{lem}\label{lem1}
Let $G$ be a primitive permutation group of degree $n$, other than
$S_n$ or $A_n$. Then either
\begin{enumerate}[(a)]
\item $G$ is $S_m$ or $A_m$ on 2-sets $(n ={m \choose 2})$, or $G$ is a
subgroup of $S_m \wr S_2$ containing $A_m^2$ ($n = m^2$); or 
\item $|G| \mbox{ \raisebox{-.9ex}{$\stackrel{\textstyle <}{\sim}$} } 
\exp (n^{1/3} \log n)$. 
\end{enumerate}
\end{lem}

We call $G$ ``large'' or ``small'' according as the first or second
alternative holds. Note that large groups have order roughly 
$\exp (n^{1/2} \log n)$.  While Lemma~\ref{lem1} uses the classification
of finite simple groups, a remarkable recent result by 
graduate students Xiaorui Sun
and John Wilmes~\cite{sun-wilmes} (extending~\cite{orders})
combined with \cite{orders2} or \cite{orders3} yields an elementary
proof of a slightly weaker bound, namely, $\exp(n^{1/3}(\log
n)^{7/3})$ in part (b), which would be just as adequate for our
purposes.

\begin{rk}
This recent progress has not entirely eliminated our dependence
on the classification of finite simple groups, and
it is worth pointing out just how the classification is used. 
First of all, we actually require stronger bounds than Lemma~\ref{lem1}
with a longer list of exceptions (\cite{fpgfsg}, see Lemma~\ref{fifth}). 
The best explicit result in this
direction is due to Mar\'oti~\cite{maroti}, but we do not need the full force
of this. We also use the classification of $2$-homogeneous groups and the facts
that simple groups can be generated by 2 elements and have small outer
automorphism groups; but these could probably be avoided with care. 
\end{rk}

	The results in this paper were mostly obtained during the Second
Japan Conference on Graph Theory and Combinatorics at Hakone in 1990. 
We are also grateful to M.~Deza, I.~Farad\v{z}ev, and M.~H.~Klin for
asking persistent questions, and especially to Klin for
proposing the question and for several contributions to the proof.
The final steps in the argument (reported in Section~8) eluded us for
some time, and so publication has been rather seriously delayed! In the
intervening time, some of the results we used in the proof have been improved
(for example in the papers ~\cite{lps,maroti}), which allows a small amount of
streamlining of our arguments; but we have kept the original arguments almost
unchanged.

\section{Outline of the proof}

	Throughout the paper we use the term ``maximal subgroup'' to
mean ``subgroup of $S_n$, maximal in the set of permutation groups
other than $S_n$ and $A_n$'': that is, a maximal subgroup of $S_n$
other than $A_n$, or a maximal subgroup
of $A_n$ contained in no other proper subgroup of $S_n$. We deduce the
Theorem from the following result.

\begin{lem}[Main Lemma]
Let $G$ be primitive on $X$ with $|X| = n$, $G \not= S_n, A_n$. Then
with probability $1- \exp (-n^{1/2 + o(1)})$, a random subset $Y$ of $X$
has the property that $M_Y = 1$ for every maximal subgroup $M$
containing $G$.
\end{lem}

(Here $M_Y$ denotes the setwise stabilizer of $Y$.)

	The deduction of Theorem~\ref{all} from the Main Lemma runs as
follows. Clearly the theorem holds for $G = S_n$; so we may assume $G
\not= S_n, A_n$. Let $Y$ be a subset for which the conclusion of the
Main Lemma holds, and let $H = \Aut(X, Y^G)$. Then $H \geq G$, and so
$H \leq M$ for some maximal subgroup $M$ containing $G$. Thus $M_Y =
1$, and so $H_Y = 1$, and
\[
|G| \geq |Y^G| = 
   |Y^H| = |H| \geq |G|,
\]
from which $H = G$ follows.

	Now to prove the Main Lemma, we need estimates for
\begin{enumerate}[(a)]
\item the number of conjugacy classes of maximal subgroups;
\item for each conjugacy class, the number of subgroups in that class
containing $G$; 
\item for each maximal subgroup $M$, the probability that $M_Y = 1$
for random~$Y$.  
\end{enumerate}

The estimate for (a) is taken from \cite{ls}, improving earlier bounds in
\cite{generating}, \cite{nmax}, and \cite{lps}. 

\begin{lem}\label{lem2}
The number of conjugacy classes of maximal subgroups is at most
$(\frac{1}{2}+o(1))n$. 
\end{lem}

	Estimates for (b) are given in the next two sections; we
separate the cases of large and small maximal subgroups. In
section~\ref{rigid} we prove a uniform bound $\exp(-c\sqrt{n})$ for
(c), using 
elementary bounds on the minimal degree.
In section~\ref{comp} we do the accounting
necessary to prove the Main Lemma. Finally in section~\ref{bound}, we
indicate how to modify the argument in order to prove 
Theorem~\ref{best}.

\medskip

We will need the following result at two points in the proof.

\begin{lem}
A primitive group of degree $n$ can be generated by at most $c(\log n)^2$
elements.
\label{gens}
\end{lem}

\begin{pf}
If $G$ has abelian socle,
then $|G| \leq n^{1+ \log n}$ and the claim is clear. So let $G$ have
socle $T^k$, where 
$T$ is non-abelian and simple.
Then $G \leq T^k \cdot
(\mbox{Out}T)^k \cdot S_k$. Now $T^k$ requires at most $2k$
generators; a subgroup of $(\Out T)^k$, at most $3k$
(since $\Out(T)$ is at most $3$-step cyclic); and a
subgroup of $S_k$, at most $2k$~\cite{chains-lb, chains-cst}. 
Since $k \leq \log n$, the result holds.
\end{pf}

\section{Large maximal subgroups}\label{large}
Let $S_m^{(2)}$ denote the action of $S_m$ on the $n={m\choose 2}$
pairs.

	In this section, we estimate the number of subgroups 
$S_m^{(2)}$
or $S_m \wr S_2$ containing a given primitive group $G$.
In the first case, we obtain a best possible result.

\begin{lem}\label{triang}
A primitive group of degree $n$ is contained in at most one subgroup
isomorphic to $S_m^{(2)}$.  
\end{lem}
~
\begin{pf}
Let $S_m$ act on $\Delta = \{1, \ldots, m \}$, and identify $X$
with the set of 2-subsets of $\Delta$. If $G \leq S_m$ on pairs and $G$
is primitive on $X$, then certainly $G$ is 2-homogeneous on $\Delta$
(transitive on unordered pairs).  These groups are listed, for
example, in \cite[pp.194--197]{hbk}.  Below we give original
sources.


\begin{itemize}
\item[(a)]
 Affine case: $G$ has an elementary abelian normal subgroup $N$ regular on
$\Delta$. Then $N$ is intransitive on $X$; so $G$ cannot be primitive.
\end{itemize}

We note that this includes the case when $G$ is not 2-transitive;
indeed, in that case $G$ has odd order and is therefore
solvable; being primitive on $\Delta$, its minimal normal
subgroup is transitive and elementary abelian.

Therefore we are left with the cases of non-affine, doubly transitive
groups $G$ on $\Delta$.  The socle $T$ of $G$ is then a non-abelian
simple group.   If $T$ is alternating of degree $k<m$ then 
$k\le 7$ and $n\le 15$ by Maillet's 19th century result~\cite{maillet};
we move these cases to item (iii) below.
The case $T=A_m$ is trivial.
In the remaining cases, $T$ is either of Lie type (these cases were
classified by Curtis, Kantor and Seitz~\cite{2trans}) or $T$ is sporadic.  
We consider each case. 
\begin{itemize}
\item[(b)]
 $G\cong\mathrm{Sp}(2d,2)$, $m = 2^{2d-1} \pm 2^{d-1}$. Then $\Delta$ is
embedded in $\mathrm{AG}(2d,2)$ in a $G$-invariant fashion, and $G$ preserves
the restriction of the parallelism relation to 2-subsets of $\Delta$; so
$G$ cannot be primitive on $X$.
\item[(c)]
 $G$ preserves a Steiner system $S(2,k,m), k>2$. (This includes groups
containing $\PSL(d,q)$ with $d>2$, unitary and Ree groups in Table~7.4,
p.197, of \cite{hbk}.) Then the set of
pairs contained in a block of the Steiner system is a block of
imprimitivity for $G$ on $X$; once again $G$ can not be primitive.
\item[(d)]
In the remaining cases, overgroups $S_m^{(2)}$ of $G$ correspond
bijectively to unions of orbital graphs of $G$ which are isomorphic
to the line graph of $K_m$. We show that, in each case, there is
at most one such union.
\begin{enumerate}[(i)]
\item $Sz(q) \leq G \leq \Aut(Sz(q))$, $m = q^2+1$.
Let $N$ be the Sylow 2-normalizer
in $Sz(q)$. Then $N$ fixes a point of $\Delta$, say 1, and is transitive
on the remaining points. Thus $N$ has an orbit of length $q^2$ on $X$,
consisting of pairs $\{1,i\}$ for $i \in \Delta \backslash \{1\}$. Any
other $N$-orbit on $X$ has length divisible by $q-1$, by considering
the 2-point stabilizer in $Sz(q)$. So the orbit of length $q^2$ is
unique. Now this orbit must be a clique in the required graph, and the
edge sets of its translates cover all edges of the graph. So the graph
is unique.
\item $\PSL(2,q) \leq G \leq \PGaL(2,q)$, $m = q+1$.
This case is similar to the preceding one.
\item Finitely many others. These are handled by \textit{ad hoc} methods.
For example, the groups $M_n$, for $n=11, 12, 23, 24$, are $4$-transitive
on $\Delta$, and so have just two non-trivial orbital graphs, the line graph
of $K_m$ and its complement.
\end{enumerate}
\end{itemize}
\end{pf}

\begin{lem}\label{lattice}
A primitive group of degree $n = m^2$ is contained in at most 
$n^{c \log n}$ subgroups isomorphic to $S_m \wr S_2$.
\end{lem}

\begin{pf}
If $G$ has an overgroup $M \cong S_m \wr S_2$, then $G$ has a
subgroup $H$ of index 2 which has a block $B$ of imprimitivity of size
$\sqrt{n}$. The pair $(H,B)$ determines $M$. (Indeed,
it is easy to show that any such pair gives rise to an overgroup of
the correct form.) So we have to estimate the number of such pairs.

\begin{enumerate}[(a)]
\item $G$ has at most $n^{c \log n}$ subgroups of index 2. This is
immediate from Lemma~\ref{gens}.
\item A block of imprimitivity containing a point $x$ for a transitive
group $H$ is determined by its stabilizer, a subgroup $K$ of $H$
containing $H_x$. Now $K$ is generated by at most $\log n$ cosets of
$H_x$, so there are at most $n^{\log n}$ choices for it.
\end{enumerate}
Multiplying (a) and (b) yields the Lemma.
\end{pf}

\section{Small maximal subgroups}\label{small}

	In this section we prove a general bound. It is good enough
for our purposes for small maximal subgroups, but not for large ones.

\begin{lem}\label{lem5}
Let $G$ be primitive, $M$ maximal, of degree $n$. Then $G$ lies in at
most $(c|M|)^{c \log^2 n}$ conjugates of $M$.
\end{lem}

\begin{pf}
First, $M$ contains at most $|M|^{c \log^2 n}$ copies of $G$,
because $G$ has at most $c \log^2 n$ generators, as we observed in 
Section 2. 

	Next, $|N_{S_n}(G) : G| \leq \exp (c \log^2 n)$. For
$|N_{S_n}(G):G|$ does not exceed $|\Out(N)|$, where $N$ is the
socle of $G$. The bound is clear if $N$ is abelian. Otherwise, $N =
T^k$, where $k \leq \log n$; and $|\mbox{Out}N| = |\mbox{Out}T|^k k$!.
If $T$ is alternating, then $|\Out T|$ is bounded by a constant;
otherwise, $|\Out T| \leq \log |T|$, and $|T| \leq n^{\log n}$.
In either case, the bound holds.

	Now let $x$ be the number of conjugates of $M$ containing $G$.
Counting pairs $(G',M')$, where $G'$ and $M'$ are conjugates of $G$
and $M$ respectively, we obtain
\[
x | S_n : N_{S_n} (G) | \leq |M|^{c \log^2 n} | S_n : M|.
\]
Rearranging, $x | M:G| \leq |M|^{c \log^2 n} |N_{S_n}(G) : G|$, which
gives the result (since $|M:G| \geq 1$).
\end{pf}

\begin{rk}
This bound is probably much too large.
\end{rk}

\section{The probability of rigidity}\label{rigid}

	The analysis in this section has been performed several times
by different people for various applications. The \textit{minimal degree}
of $G$ is the least number of points moved by a non-identity element of
$G$.  
Doubly transitive groups other than $S_n$ and $A_n$ have minimal
degree at least $n/4$ by Alfred Bochert's 1892 combinatorial 
gem~\cite{bochert}.  Primitive but not
doubly transitive groups have minimal degree at least $(\sqrt{n}-1)/2$
by elementary arguments~\cite[Thm 0.3]{orders}, so the same 
lower bound\footnote{Still by elementary arguments, the lower
bound $2\sqrt{n}$ holds for the minimal degree of all primitive groups
other than $S_n$ and $A_n$, 
for all sufficiently large $n$~\cite{sun-wilmes}.}
holds for all primitive groups other than $S_n$ and $A_n$.  
For convenience we cite the following
slightly stronger bound taken from~\cite{mindeg}.

\begin{lem}\label{lem6}
A primitive permutation group of degree $n$, 
other than $S_n$ or $A_n$, 
has minimal degree at least $(\sqrt{n})/2$.
\end{lem}

\begin{lem}\label{lem7}
Let $M$ be a maximal primitive group of degree $n$, acting on a set
$X$. If $Y$ is a random subset of $X$, then 
\[
\mbox{Prob}(M_Y \not= 1) \leq \exp (-c \sqrt{n}),
\]
for some constant $c$.
\end{lem}

\begin{pf}
Again we treat large and small groups separately. If $M$ is
large, we require the probability that a random graph (or bipartite
graph) admits a non-trivial automorphism, for which estimates exist
\cite{asymm}.

	So suppose $M$ is small. Let $m$ be its minimal degree. If
$g \in M$, $g \not= 1$ then $g$ has at most $n - m/2$ cycles
on $X$ (the extreme case occurring if $g$ is an involution moving $m$
points), and so $g$ fixes at most $2^{n - m/2}$ subsets. So the
probability that a random subset is fixed by $g$ is at most $2^{-m/2}
\leq 2^{-\sqrt{n}/4}$. Then
\begin{eqnarray*}
\Prob(M_Y \not= 1) & \leq & |M| \cdot 2^{-\sqrt{n}/4} \\
			& \leq & \exp (n^{1/3} \log n) \cdot 2^{-\sqrt{n}/4} \\
			& \leq & \exp (-c \sqrt{n}).
\end{eqnarray*}
\end{pf}

\section{Completion of the Proof}\label{comp}

	Now by the above Lemmas, the number $F(G)$ of maximal subgroups
containing $G$ is at most
\[
1+n^{c \log n} + \exp (\log^4 n) \cdot (c \exp (n^{1/3} \log n))^{c
\log^2 n},
\]
and the probability that the conclusion of the Main Lemma fails is at
most $F(G) \exp (-c\sqrt{n})$. So the result is proved.

\section{Bounding the size of sets required}\label{bound}

	The proof of Theorem~\ref{best} follows closely the argument we
have given, but the technical details are considerably harder.
The difficulty arises because the analogue of
Lemma~\ref{lem7} is much weaker. If $M$ is maximal primitive with
minimal degree $m$,   
then a non-identity element of $M$ has at most
$n-m/2$ cycles, and so fixes at most
\[\sum_{i=0}^k{n- m/2\choose i}\leq 2{n-m/2\choose k}\]
subsets of size $k$. 
So
\begin{eqnarray*}
\Prob(M_Y \not= 1) & \leq & 2|M|  {n- m/2 \choose k} \mbox{\Large $/$}
{n \choose k}\\ 
& \leq & 2|M| (1-\frac{m}{2n})^k\\
& <    & 2|M| \exp (-\frac{km}{2n}).  
\end{eqnarray*}

	Now $|M| \geq 2^{n/m}$ 
(this bound holds for any transitive group of degree $n$ and minimal 
degree $m$~\cite{mindeg, mindeg2}), 
and so
\begin{eqnarray*}
\Prob(M_Y \not= 1) & \leq & |M| \exp\left(-\frac{k}{2 \log |M|}\right)\\
& = & 2\exp\left(\log|M|-\frac{ck}{\log|M|}\right).
\end{eqnarray*}
For $k=n^{1/2+\epsilon}$, no conclusion is possible unless 
$|M| \leq \exp (n^{1/4- \epsilon})$.

	Fortunately the classification of finite simple groups gives
such a bound with known exceptions. For this result see \cite{fpgfsg,maroti},
but note that we do not need the full refinement of Mar\'oti's estimates.

\begin{lem}\label{fifth}
If $G$ is primitive of degree $n$ and maximal, then either
\begin{enumerate}[(a)]
\item $G$ is contained in $S_m$ on $k$-sets, $k= 2$, $3$ or $4$, ($n =
{m \choose k}$), or $S_m \wr S_k$, $k=2$, $3$, or $4$ ($n = m^k$); or
\item $|G| \leq \exp(n^{1/5} \log n)$.  
\end{enumerate}
\end{lem}

The groups contained in $S_m\wr S_2$, where $S_m$ acts on $2$-sets
(with order around 
$\exp(n^{1/4}\log n)$) 
do not need to  be considered since
they are contained in $S_{m(m-1)/2}\wr S_2$.

	If we redefine ``large'' maximal subgroups to include all
those under (a), then our estimates suffice for ``small'' maximal
subgroups, and we can use separate estimates for the probabilities
that random $k$-uniform hypergraphs and random $k$-partite $k$-uniform
hypergraphs, for $k=2, 3, 4$, admit non-trivial automorphisms.
These are given in the next Section. We also need to prove analogues of 
Lemmas~\ref{triang} and~\ref{lattice} for $k=3,4$.

\begin{lem}
A primitive group of degree $n$ lies in at most one copy of $S_m$ on
$k$-sets, for fixed $k\ge2$ and $m>2k$.
\end{lem}

\begin{pf}
The case $k=2$ is Lemma~\ref{triang}. For $k\ge3$, there are
very few groups other than $S_m$ and $A_m$ which act primitively
on $k$-sets. For such a group must be $k$-homogeneous, and hence
$2$-transitive~\cite{lw}; as before, it cannot have a regular normal subgroup
$N$ (since $N$ would be intransitive on $k$-sets). This leaves
only the cases $\PSL(2,q)\le G\le\PGaL(2,q)$ (with
$k=3$) or $G$ is a Mathieu group. In the first case, the stabiliser
of a $3$-set is not maximal (by inspection of the list of maximal
subgroups in~\cite{dickson}). The Mathieu groups are handled by
\textit{ad hoc} methods.
\end{pf}

\begin{lem}
A primitive group of degree $n$ lies in at most $n^{c\log n}$
copies of $S_m\wr S_k$, for fixed $k\ge 2$ and $m>2$.
\end{lem}

\begin{pf}
The proof of Lemma~\ref{lattice} applies with trivial changes.
\end{pf}

\section{Asymmetry of random hypergraphs}

\begin{lem}
For $2\le t\le n/2$, let $P_t$ denote the probability 
that a random $t$-uniform hypergraph on $n$ vertices
is not asymmetric. For $t=2$ we have $P_2 = \sqrt{2}n^22^{-n/2}(1+o(1))$.
For $t\ge 3$ we have $P_t < \exp(-c_1{n-1 \choose t-1}) < 2^{-cn^2}$
for some positive absolute constants $c,c_1$.
\end{lem}

\begin{pf}
(For $t=2$ this is well known \cite{asymm}, but we include this
case, too, in the proof.) If the random hypergraph $\mathcal{H}$ has a 
nonidentity automorphism then it has one of prime order. Let $\sigma$ be a
permutation of $V(\mathcal{H})$ of prime order $p$. 
Let $N(\sigma)$ denote the number of $t$-sets moved by $\sigma$, so that
$N(\sigma)\le {n\choose t}$. The number of cycles of $\sigma$ on $t$-sets
is at most ${n\choose t}-N(\sigma)/2$, the extremal case being where $\sigma$
is an involution.

Let $P(\sigma)$ denote the probability that $\sigma$ is an automorphism
of $\mathcal{H}$. We think of $\mathcal{H}$ being generated by flipping a
coin for every $t$-tuple 
$F\subset V$ to decide whether or not to include $F$ in $E(\mathcal{H})$.
However, for each $\sigma$-orbit of $t$-sets, we can flip the coin
only once. Therefore $P(\sigma)\le 2^{-N(\sigma)/2}$ and
\begin{equation} \label{est.eq}
  P_t \le \sum_{\sigma} 2^{-N(\sigma)/2},
\end{equation}
where the summation extends over all permutations $\sigma$ of
prime order.

Let $s$ denote the size of the support of $\sigma$ (number of elements
of $V(\mathcal{H})$ moved). Note that $p|s$.
Let $\rho(s,\ell,p)={s/p \choose \ell/p}$ if $p$ divides $\ell$, and $0$ 
otherwise. Let us compute $N(\sigma)$ by counting for each $\ell\ge 1$
those $t$-sets which intersect the support of $\sigma$ in exactly
$\ell$ elements. Adding these up we obtain
$$ N(\sigma) =\sum_{\ell =1}^t \left({s\choose \ell}-\rho(s,\ell,p)\right)
    {n-s\choose t-\ell}.$$
This quantity is estimated as
\begin{equation}  \label{est2.eq}
 N(\sigma) \ge {1\over 2}\sum_{\ell=1}^t {s\choose \ell}{n-s\choose t-\ell}
  = {1\over 2} \left({n\choose t}-{n-s\choose t}\right).
\end{equation}
For $t=2$, the right hand side is 
$$ N(\sigma) \ge \frac{1}{2}\left({n\choose 2} - {n-s\choose 2}\right) 
= s(n-s/2-1/2)/2.$$
>From equation (\ref{est.eq}) we then infer that 
$$ P_2 \le \sum_{s=2}^n {n\choose s}(s!-1)2^{-s(n-s/2-1/2)/4}=
   {n\choose 2}2^{(-n+3/2)/2}(1+o(1)),$$
as stated. The opposite inequality follows by taking the second
term of the inclusion-exclusion formula into account in calculating
the probability that there exist two vertices switched by
a transposition.

For $t\ge 3$, equation (\ref{est2.eq}) implies $N(\sigma)>{n-1\choose t-1}
> c_2 n^2$, therefore
$$ P_t \le n! 2^{-(1/2){n-1\choose t-1}}.$$
\end{pf}

A $t$-partite \textit{transversal hypergraph}\/
$\mathcal{H}$ is a $t$-uniform hypergraph whose
vertex set is partitioned into $t$ 
\textit{``layers''} $V(\mathcal{H})=V_1\cup\cdots\cup V_t$
and each edge intersects each class $V_i$ in exactly one element. The partition
into layers is given and is part of the definition of $\mathcal{H}$.
Automorphisms of $\mathcal{H}$ preserve the partition by definition (but may
interchange the layers). A transversal hypergraph is \textit{balanced}\/
if all layers have equal size.

\begin{lem}
For $2\le t\le n/2$, let $Q_t$ denote the probability 
that a random balanced $t$-partite transversal hypergraph on $n=tr$ vertices
is not asymmetric. For $t=2$ we have $Q_2 = cn^22^{-n/4}(1+o(1))$.
For $t\ge 3$ we have $Q_t < \exp(-c_1{n-1 \choose t-1}) < 2^{-cn^2}$
for some positive absolute constants $c,c_1$.
\end{lem}

\begin{pf} 
As before, let $\sigma$ be a permutation of $V=V(\mathcal{H})$ of prime order 
$p$. By our remark about automorphisms before the lemma, $\sigma$ respects the
partition into layers $(V_1,\ldots,V_t)$.
We say that a $t$-tuple is \textit{transversal} if it intersects
each layer in exactly one element. The total number of transversal
$t$-sets is $r^t$ where $r=n/t$.
Let $N(\sigma)$ denote the number of transversal $t$-sets
moved by $\sigma$. Let $Q(\sigma)$ denote the probability that
$\sigma$ is an automorphism of the random transversal hypergraph $\mathcal{H}$.
Just as in equation~(\ref{est.eq}), we have
\begin{equation} \label{estt.eq}
  Q_t \le \sum_{\sigma} 2^{-N(\sigma)/2},
\end{equation}
where the summation extends over all permutations $\sigma$ of
prime order, respecting the partition $(V_1,\ldots,V_t)$.

If $\sigma$ moves some of the layers then it moves at least $p$ 
layers.  
Hence in this case,
\begin{equation}  \label{estt1.eq}
N(\sigma) \ge r^t - r^{t-p+1} > r^t/2.
\end{equation}
If $\sigma$ fixes all layers, then let $s_i$ denote the number of
elements in $V_i$ that are moved by $\sigma$. So the support of $\sigma$
has size $s=s_1+\cdots+s_t$.

The transversal $t$-sets not moved by $\sigma$ are now exactly those
which do not intersect the support of $\sigma$. Therefore
\begin{equation}  \label{estt2.eq}
 N(\sigma) = r^t - \prod_{i=1}^t (r-s_i) \ge r^t(1-\eee^{-s/r})
   \ge r^t\cdot {s\over r}\cdot \left(1-{s\over 2r}\right).
\end{equation}
We infer that for $s \le r/2$ we have
\begin{equation}
 N(\sigma) \ge r^t 3s/(4r) \ge r^{t-1},
\end{equation}
and for $s\ge r/2$ we have (for $r\ge 3$)
\begin{equation}
 N(\sigma) \ge r^t (1-\eee^{-1/2}) > 0.39\cdot r^t > r^{t-1}.
\end{equation}
For $t=2$ we conclude, separating the case $s=2$, that 
$$ Q_2 \le  (({n\over 2})!)^2 2^{-n^2/8} + {n(n-2)\over 4} 2^{-n/4}
    + \sum_{s=4}^{n/4} {n\choose s} 2^{-3ns/16}+
      \sum_{s=n/4}^n {n\choose s} 2^{-0.39n^2/2}.$$
It is clear that the second term dominates this sum.

For $t\ge 3$ we obtain $Q_t \le n! 2^{-r^{t-1}/2}.$
\end{pf}

\section{Open problems}  
\medskip\noindent
1.  Can we eliminate the need for the classification
of finite simple groups from the proof of Theorem~\ref{best}?

\bigskip\noindent
2.  Assuming $G$ is a primitive group whose socle is not 
a product of alternating groups, is it possible to reduce
the size of the hyperedges in Theorem~\ref{best} below
$n^{0.49}$ (for sufficiently large $n$)?   Is it possible
to reduce it to $n^{o(1)}$ ?

\end{document}